\documentstyle{amsppt}
\magnification=1200
\hsize=150truemm
\vsize=224.4truemm
\hoffset=4.8truemm
\voffset=12truemm

\NoRunningHeads

\define\C{{\bold C}}
\let\th\proclaim
\let\fth\endproclaim

\newcount\refno
\global\refno=0

\def\nextref#1{
      \global\advance\refno by 1
      \xdef#1{\the\refno}}

\def\bref {\ref\global\advance\refno by 1\key{\the\refno}}

\nextref\AL
\nextref\BED
\nextref\BK
\nextref\CHI
\nextref\DUV
\nextref\DS
\nextref\ELI
\nextref\FOR
\nextref\FORS
\nextref\GRO
\nextref\IS
\nextref\SLO
\nextref\STO
\nextref\WHI
\nextref\YE

\topmatter
\title 
Polynomial hull of a torus fibered over the circle
 \endtitle
\author Julien Duval \footnote{Laboratoire de Math\'ematiques, Univ.Paris-Sud, Universit\'e Paris-Saclay, Orsay, France \newline julien.duval\@math.u-psud.fr\newline}
 and Mark Lawrence 
\footnote{Department of Mathematics, Nazarbayev University, Astana, Kazakhstan \newline mlawrence\@nu.edu.kz\newline}
\footnote""{Keywords: totally real torus, polynomial hull, holomorphic disc \newline AMSC : 32E20 \newline}
\endauthor
\abstract Given a 2-sheeted torus over the circle with winding number $1$, we prove that its polynomial hull is a union of 2-sheeted holomorphic discs. Moreover when the hull is non degenerate its boundary is a Levi-flat solid torus foliated by such discs.
\endabstract 
  \endtopmatter 
\document
 \subhead 0. Introduction \endsubhead

\null

Let $M$ be a totally real torus in $\C^2$. We are interested in the structure of its polynomial hull $\hat M$. We specialize in the situation where $M$ sits in the vertical cylinder $S^1\times \C$ over the unit circle. We suppose moreover that the vertical projection $\pi$ restricts to a submersion from $M$ to $S^1$. Its fibers $M_z$ are $n$ curves moving around smoothly when $z$ winds around $S^1$. For short we speak of a {\it $n$-sheeted torus}. In $\hat M$ there is a trivial part $Y=\hat M\vert_{S^1}$ over the circle. It consists in filling in the curves of $M_z$. So $Y$ is a solid torus (diffeomorphic to $S^1\times D^2$) fibered over the circle with fibers $Y_z$ consisting in $n$ discs. Call it the {\it vertical filling} of $M$. Note that $\hat M=\hat Y$. We want to describe the rest of $\hat Y$, the part $\hat Y\vert_D$ over the open unit disc $D$. 

\null
When the torus is $1$-sheeted the situation is well understood thanks to works by Forstneric [\FORS] and Slodkowski [\SLO]. Namely $\hat Y\vert _D$ (when non empty) consists in a union of holomorphic discs with boundary in $M$ which are graphs over $D$. For short we speak of {\it 1-sheeted discs}. If moreover $\hat Y\vert _D$ is {\it non degenerate} (not reduced to a single disc) then its boundary $\partial \hat Y\vert_D$ is a Levi-flat solid torus foliated by such discs, bounding $M$. We call it the {\it (horizontal) filling} of $M$.

\null
In this paper we take up the next simplest case, namely when the torus is $2$-sheeted. Define its {\it winding number} $w(M)$ as follows. Take a section of $Y$, a continuous map $\gamma : S^1\to Y$ of the form $\gamma(z)=(z^2,\alpha(z))$. Then $w(M)$ (necessarily odd) is the winding number of $ \alpha(z)-\alpha(-z)$ around $0$. 
    Here is our result which parallels the 1-sheeted case (see also [\WHI] for a particular case).

\th {Theorem } Let $M$ be a 2-sheeted torus with $w(M)=1$. Then $\hat Y\vert_D$ (when non empty) consists in holomorphic discs with boundary in $M$ which are 2-sheeted over $D$. When $\hat Y\vert_D$ is non degenerate (not reduced to a single disc) then its boundary $\partial \hat Y\vert_D$ is a Levi-flat solid torus foliated by such discs, bounded by $M$. \fth

Here a {\it 2-sheeted disc} is a holomorphic map $f: D\to \C^2$ such that $\pi \circ f:D \to D$ is a Blaschke product of degree $2$. Moreover $f$ extends smoothly up to $\partial D$ and $f(\partial D) \subset M$. We again call the solid torus obtained in the non degenerate case the (horizontal) filling of $M$.
   
\null
The proof of the theorem goes by a continuity method as for the 1-sheeted case. The description of the boundary of the hull requires however a new ingredient, a surgery process. The rough idea is to start with a big 1-sheeted torus and to shrink it down to $M$ passing through a singular torus. According to [\SLO] the initial torus bounds a filling which deforms up to the singular level. We then get graphs with boundary into the singular torus. The point is that they are touching once by pair at the boundary. We now make a surgery on the pairs to get a family of 2-sheeted discs over $D$. They build a solid torus whose boundary is a new 2-sheeted torus. Shrinking this torus down to $M$ the filling again follows, giving at the limit either the filling of $M$ or the degeneracy of the hull.  

\null

Our guiding model is $(\vert w^2-z\vert=t)$ over $S^1$. For $t>1$ we have a 1-sheeted torus, which becomes singular for $t=1$ and 2-sheeted (with winding number 1) for $t<1$. Correspondingly for $t>1$ the filling $(\vert w^2-z\vert=t)$ over $D$ is foliated by graphs $(w=\pm \sqrt {z+tu})$ ($u$ in $S^1$), which touch by pair at the boundary for $t=1$, and merge into 2-sheeted discs $(w^2=z+tu)$ for $t<1$.

\null
Before entering the details we collect some background. 

\subhead 1. Background \endsubhead

\null
In this section $M$ stands for a torus over $S^1$, either 1- or 2-sheeted, and $Y$ for its vertical filling. We also use deformations $(M_t)$ of $M$ among tori of the same type. Everything is smooth except otherwise mentioned. 

\null\noindent
a) {\bf Polynomial hull.}

The polynomial hull $\hat K$ of a compact set $K$ in $\C^2$ is the set of points $x$ such that $\vert P(x)\vert \leq \Vert P\Vert_K$ for all polynomials $P$ (see [\STO]). The set $K$ is polynomially convex if $\hat K=K$. In this case $K$ satisfies Runge theorem. Note that any holomorphic disc with boundary in $K$ (for short {\it attached to $K$}) is contained in $\hat K$ by the maximum principle. The fact that $\hat K$ can be entirely described by such discs is far from true in general, but it will be the case in our situation.

\null\noindent
{\bf Oka principle.} 
 Geometrically $\C^2\setminus\hat K$ is a union of algebraic curves which can be swept out continuously toward infinity. This is the content of Oka principle (see [\STO]). It can be localized and translates in our situation as follows. 

Consider a continuous family of holomorphic discs (or pairs of discs) $(\Delta_t)_{t\in[0,1]}$ over $D$ such that $\partial \Delta_t\subset S^1\times \C$ avoids $M$ for all $t$ and $\Delta_0$ avoids $ \hat Y$. Then $\Delta_t$ avoids $\hat Y$ for all $t$.

\null\noindent
{\bf Rational convexity.} 
Replacing polynomials by rational functions we may also speak of rational hull and rational convexity.
It turns out that in our situation the torus $M$ is rationally convex (see [\DS]). As it also is totally real (never tangent to a complex line) it actually becomes lagrangian for a K\"ahler form $\omega$ in $\C^2$ (see [\DUV]), meaning that $\omega\vert_M$ vanishes.

\null\noindent
b) {\bf Holomorphic discs.} 

We investigate the holomorphic discs attached to $M$. They automatically are smooth up to the boundary (see [\CHI]). We are mostly interested in embedded discs which are 1- or 2-sheeted over $D$ accordingly to $M$. For short we call them {\it horizontal discs}.

\null\noindent
{\bf Area.} As $M$ is lagrangian for $\omega$ we get a bound for the area of a holomorphic disc $\Delta$ attached to $M$. Indeed consider a primitive $\lambda$ of $\omega$. Then area($\Delta)\leq C\int_\Delta \omega =C\int_{\partial\Delta}\lambda$. As $\lambda\vert_M$ is closed this bound only depends on the homology class of $\partial \Delta$ in $M$. This still holds along a deformation of such tori.

\null\noindent
{\bf Index.} 
Given an oriented loop $\gamma$ in $M$ we define its index $\mu(\gamma)$ as follows (see [\FOR],[\AL], our index is half of Maslov index). Consider a continuous frame $(t(x), n(x))$ of $T_xM$ when $x$ goes around $\gamma$. Then $\mu(\gamma)$ is the winding number of $\det(t,n)$ around $0$. This makes sense as $M$ is totally real. It only depends on the homology class of $\gamma$. Conversely this class is completely determined by $\mu(\gamma)$ and the degree of $\pi \circ \gamma$. For a holomorphic disc $\Delta$ attached to $M$ define $\mu(\Delta)$ as $\mu(\partial \Delta)$. 

When $\Delta $ is horizontal $\mu(\Delta)$ can also be computed geometrically as follows. Take a continuous section $\gamma$ of the interior of $Y$. Then $\mu(\Delta)-1$ is the linking number of $\partial \Delta$ and $\gamma$. In other words if $\gamma$ bounds a disc (not necessarily holomorphic) $\Delta'$ over $D$ then $\mu(\Delta)-1$ coincides with the intersection number $\Delta.\Delta'$. 

\null\noindent
{\bf Deformation.} 
This index $\mu$ gives the size of the possible deformations of a horizontal disc $\Delta$ attached to $M$ (see [\FOR],[\AL]). Namely nearby horizontal discs attached to $M$ form a manifold of dimension $2\mu-1$ when $\mu\geq 1$. Here we consider the discs geometrically, moding out their reparametrizations. When $\mu=1$ these discs are actually disjoint. On the other hand when $\mu>1$ we are able to deform $\Delta$ keeping a prescribed fixed point.

Similar statements hold when we are given a deformation $(M_t)$ of $M=M_0$. Nearby horizontal discs attached to $M_t$ for small $t$ still form a manifold of dimension $2\mu-1$, smoothly depending on $t$. When $\mu>1$ we are able to deform $\Delta$ in nearby horizontal discs $\Delta_t$ attached to $M_t$ for small $t$ keeping a prescribed fixed point. 

\null\noindent
{\bf Uniqueness.} When $\mu(\Delta)<1$ the disc $\Delta$ is isolated among the horizontal discs attached to $M$.
It turns out that in this case there is no other horizontal disc attached to $M$ or even to $Y$, so $\Delta$ is unique.

If not consider another horizontal disc $\Delta'$ with $\partial \Delta'\subset Y$. Suppose first $\partial \Delta'$ contained in the interior of $Y$. Then $\Delta.\Delta'=\mu(\Delta)-1<0$ which is impossible as holomorphic discs intersect positively. In the general case $\partial \Delta'\cap \partial \Delta$ is finite. Now slightly perturb $\Delta'$ in a new holomorphic disc $\Delta''$ so that $\partial \Delta''$ enters the interior of $Y$ near the points of $\partial \Delta'\cap \partial \Delta$ (see [\SLO],[\FORS] for similar arguments). Then $\partial \Delta''$ is homotopic to a section of the interior of $Y$ in $S^1\times \C\setminus \partial \Delta$. So again $\Delta.\Delta''=\mu(\Delta)-1<0$ which is impossible.

\null\noindent
{\bf Compactness.} Consider a sequence of horizontal discs $\Delta_n$ attached to $M$ of fixed index $\mu$. Then the boundaries $\partial \Delta_n$ are in a fixed homology class of $M$. Hence the discs $\Delta_n$ are of bounded area. By Gromov compactness theorem ([\GRO],[\IS]) after extracting $\Delta_n$ converges to a bunch of holomorphic discs $\Delta_\infty \cup_i \Delta'_i$ attached to $M$. Here $\Delta_\infty$ is horizontal and the discs $\Delta'_i$ are vertical (they project down to points $z_i$ in $S^1$). We call them the {\it vertical bubbles}. 

Moreover $\partial \Delta_\infty\cup_i \partial \Delta'_i$ is homologous to $\partial \Delta_n$. So $\mu=\mu_\infty+\sum_i \mu'_i$ for the indices. Note that $\mu'_i\geq 1$ as it is the winding number of the tangent of $M_{z_i}$ (or one of its components in the 2-sheeted case). Therefore $\mu_\infty\leq\mu$ with strict inequality in the presence of vertical bubbles. When $\mu=1$ bubbles don't appear. If they did we would get a horizontal limit of index $<1$ contradicting the presence of other discs attached to $M$. In this case the convergence is nice (smooth in adequate parametrizations).

This compactness also holds when we have a deformation $(M_t)_{t\in [0,1]}$ of $M=M_0$. Given a family of horizontal discs $\Delta_t$ attached to $M_t$ for $t<1$ and of fixed index we get a limit $\Delta_1\cup_i \Delta'_i$ attached to $M_1$ with the same features.

\null\noindent
{\bf Persistence of intersection.} The first form is classical. Let $\Delta_n$ be a sequence of holomorphic discs converging to a holomorphic disc $\Delta$ (or a bunch of holomorphic discs). Suppose $\Delta$ has an interior intersection with an algebraic curve $C$. Then $\Delta_n\cap C\neq \emptyset$ for large $n$.  

The second form takes into account the boundary intersections. Let $\Delta,\Delta'$ be distinct horizontal discs attached to $M$ such that $\Delta\cap \Delta'\neq \emptyset$. Suppose that $\Delta= \Delta_0$ deforms in a continuous family of discs $\Delta_t$ attached to $M$. Then $\Delta_t \cap \Delta'\neq \emptyset$ for $t$ small. 

This holds true even if the only intersections between $\Delta$ and $\Delta'$ are contacts at the boundary. They create interior intersections for small $t>0$. Actually this can be quantified. We may define an interior intersection number $i$ and a boundary intersection number $b$ between $\Delta$ and $\Delta'$. Then $b+2i$ is deformation invariant (see [\YE]).

\null\noindent
c) {\bf Filling.} 

A {\it (horizontal) filling} of $M$ is a solid torus $\Sigma$ bounded by $M$ and foliated by horizontal discs. Note that they necessarily are of index 1. 

\null\noindent
{\bf Generation.}  Conversely given a horizontal disc $\Delta$ with $\mu(\Delta)=1$ then it generates a filling. 

Indeed use local deformations and compactness to build a smooth family of horizontal discs of index 1 $(\Delta_t)_{t\in [0,1]}$, starting with $\Delta_0=\Delta$, developing on one side of $\partial \Delta$ on $M$ at the boundary and ending with $\partial \Delta_1$ touching $\partial \Delta$ from the other side for the first time. We claim that $\Delta_1=\Delta_0$ so we get the filling. If not by persistence of intersections the contact would propagate back to $\Delta_t \cap \Delta_0\neq \emptyset$ for small $t$, contradicting the disjointness of local deformations of discs of index 1 (see [\ELI],[\BK] for similar arguments).

\null\noindent
{\bf Deformation.} Given a deformation $(M_t)$ of $M=M_0$ and a filling $\Sigma_0$, then $\Sigma_0$ deforms in a filling $\Sigma_t$ of $M_t$ for small $t$ (see [\BED]). Moreover if the $M_t$ foliate a neighborhood of $M$ in $S^1\times \C$ then the $\Sigma_t$ also foliate a neighborhood of $\Sigma_0$ in $D\times \C$.  

Let now $(M_t)_{t\in [0,1]}$ be a {\it decreasing deformation}, meaning that the $Y_t$ decrease and the $M_t$ foliate $Y_0\setminus Y_1$. Suppose we have a corresponding deformation of fillings $(\Sigma_t)_{t\in[0,1]}$ such that $\Sigma_0=\partial \hat Y_0\vert _D$. Then $\Sigma_t=\partial \hat Y_t\vert _D$ for each $t$.

 Indeed on one hand $\Sigma_t \subset \hat Y_t$ by the maximum principle. On the other hand the $\Sigma_{t'}$ foliate a 1-sided neighborhood of $\Sigma_t$ for $t'<t$ and are outside $\hat Y_t$ by Oka principle.

\null\noindent
{\bf Uniqueness.} When $\partial \hat Y\vert_D$ is a filling $\Sigma$ of $M$ then it is unique. Actually there is no other horizontal disc of index 1. 

If not consider such a disc $\Delta$ distinct from the discs $\Delta_t$ foliating $\Sigma$. Note that $\partial \Delta$ and the $\partial \Delta_t$ are homologous, as the discs are horizontal of the same index. As the $\partial \Delta_t$ foliate $M$ this forces a contact between $\partial \Delta$ and one $\partial \Delta_t$. We then get an interior intersection between $\Delta$ and a nearby $\Delta_t'$ by persistence of intersection. In other words $\Delta$ crosses $\Sigma$ over $D$. But this impossible as $\Delta\subset \hat Y$ and $\Sigma=\partial\hat Y\vert_D$.

\null\noindent
{\bf Compactness.} Consider a decreasing deformation $(M_t)_{t\in[0,1]}$ with a corresponding deformation of fillings $(\Sigma_t)$ for $t<1$. Suppose that $\Sigma_0=\partial \hat Y_0\vert_D$. Then either $\Sigma_t$ converges to a filling $\Sigma_1$ of $M_1$ such that $\Sigma_1=\partial \hat Y_1\vert_D$, or $\hat Y_1\vert_D$ is degenerate, it reduces to a single horizontal disc. 

Indeed take a sequence of discs $\Delta_{t_n}\subset \Sigma_{t_n}$ when $t_n\to 1$. After extracting it converges toward a bunch of discs $\Delta\cup_i \Delta'_i$, where $\Delta$ is horizontal and the $\Delta'_i$ are vertical bubbles. If bubbles do appear then $\mu(\Delta)<1$ so $\Delta$ is the only horizontal disc attached to $M_1$. There is no other possible horizontal limit of such sequences, so $\Sigma_t=\partial \hat Y_t\vert_D$ shrinks down toward $\Delta$. Hence $\hat Y_1\vert_D =\Delta$ and the hull is degenerate. When there is no bubble in the limit then $\mu(\Delta)=1$. It generates a filling $\Sigma_1$ which can be deformed back in fillings $\Sigma'_t$ of $M_t$ for $t$ close to $1$. By uniqueness $\Sigma'_t$ coincides with $\Sigma_t$. So $\Sigma_t$ converges to $\Sigma_1$. 

\null\noindent
{\bf The 1-sheeted case.} We end up by sketching the main part of Slodkowski theorem (see [\SLO]). Let $M$ be 1-sheeted with $\hat Y\vert_D\neq \emptyset$. Then either $\partial \hat Y\vert_D$ is a filling or $\hat Y\vert_D$ is degenerate.

Indeed construct a decreasing deformation $(M_t)_{t\in [0,1]}$ from $M_0=S^1\times RS^1$ ($R$ large) to $M_1=M$. Note that $M_0$ has a trivial filling $\Sigma_0=D\times RS^1$ which coincides with $\partial \hat Y_0\vert_D$. Also note that $\hat Y_t\vert_D$ is non degenerate for $t<1$ as $Y$ sits in the interior of $Y_t$. Hence using  local deformations and compactness we construct a global deformation of fillings $(\Sigma_t)$ of $M_t$ for $t<1$. Passing to the limit we get the alternative.

\null\noindent

\subhead 2. Proof of the theorem\endsubhead

\null

From now on $M$ stands for a 2-sheeted torus with $w(M)=1$, $Y$ is its vertical filling and we assume $\hat Y\vert_D\neq \emptyset$. 

We start by the part concerning the boundary of the hull.

\null\noindent
a) {\bf Boundary of the hull.}

We prove that either $\partial \hat Y\vert_D$ is a filling $\Sigma$ of $M$ or $\hat Y\vert_D$ is degenerate.

We first sketch the argument. In what follows the geometric constructions take place over $S^1$. We modify the model of the introduction to get a smoothly immersed singular torus.

\null\noindent
{\bf Sketch.}
 The first step is to construct an immersed 1-sheeted torus $N$ out of $M$. For this we glue to $Y$ a M\"obius strip $B$. Then we smooth out the corners to get a set $Z$ and put $N=\partial Z$. It is a 1-sheeted torus pinched along $B$. By [\SLO] we have a filling of $N$. It turns out that its holomorphic discs intersect by pairs along the core $\gamma$ of $B$.

 The second step is to make a surgery on a pair of such discs to create a 2-sheeted holomorphic disc $\Delta$. Then we insert $\partial \Delta$ into a new 2-sheeted torus $M'$ obtained by opening up $Z$ along a neighborhood of $\gamma$, thickening what remains and taking the boundary. It turns out that $\mu(\Delta)=1$ so $\Delta$ generates a filling of $M'$.

 The third step is to deform back $M'$ to $M$. This builds a deformation of fillings which ends up with the desired filling $\Sigma$ or the degeneracy of the hull.

Before entering the details we need a reduction.
 
\null\noindent
{\bf Reduction.}
The situation reduces to the case where $Y$ winds once around the constant disc $\Delta_0=(w=0)$. In other words $Y$ avoids $\Delta_0$ and if $(z^2,\alpha(z))$ is a continuous section of $Y$ then the winding number of $\alpha(z)$ around $0$ is $1$.

 Actually it is enough to exhibit a holomorphic disc $(z,f(z))$ around which $Y$ winds once. We may then apply the automorphism $(z,w)\mapsto (z, w-f(z))$ of $\overline{D}\times \C$ to get the reduction. This disc is obtained thanks to Slodkowski theorem [\SLO]. 

For this consider a 2-sheeted disc $\Delta=(w^2=Rz)$ with $R$ large so that $\partial \Delta$ avoids $Y$ by far. Glue to $\partial \Delta$ a M\" obius strip avoiding $Y$. So in each fiber we connect the 2 points of $\partial \Delta$ by an arc avoiding the 2 components of $Y_z$. Moreover if we take continuous sections $(z^2,\alpha(z))$ of $Y$ and $(z,\beta(z))$ of the M\" obius strip we want the winding number of $\alpha(z)-\beta(z^2)$ around $0$ to be $1$. This can be done because $w(M)=1$. Now thicken slightly the M\" obius strip and take the boundary. We get a 1-sheeted torus with a non degenerate hull. Therefore we have plenty of 1-sheeted discs with boundary in it (see [\SLO] or \S 1 c)). By construction $Y$ winds once around them.

\null\noindent
{\bf Constructing $N$.} From now on we suppose that $Y$ winds once around $\Delta_0$. After dilation we may suppose $Y\subset (\vert w\vert >\sqrt2)$. We glue to $Y$ a M\" obius strip $B$ which is standard in $(\vert w\vert \leq\sqrt2)$. There $B=(w^2=tz)$, $(z\in S^1, t \in [0,2])$. So in each fiber we connect the 2 components of $Y_z$ by an arc $B_z$. Then we smooth the corners of $Y\cup B$ by thichkening it sligthly except in $(\vert w\vert \leq 1)$. Call the resulting set $Z$. In each fiber $Z_z$ looks like a smooth dumbbell. Take $N=\partial Z$. It is a smooth 1-sheeted torus pinched along the standard M\" obius strip in $(\vert w \vert \leq 1)$. Denote by $\gamma=\partial \Delta_0$ the core of $B$. Note that the algebraic curve $C=(w^2=2z)$ intersects $Z$ in its interior.

\null\noindent
{\bf Filling.} We verify that $N$ has a filling. 

For this consider a decreasing deformation $(N_t)_{t\in [0,1]}$ from a big standard torus $N_0=S^1\times R S^1$ to $N_1=N$. According to [\SLO] $N_t$ has a filling $\Theta_t=\partial \hat N_t\vert_D$ for $t<1$. We claim that $\Theta_t$ smoothly converges when $t\to 1$. Actually if we take a sequence of holomorphic discs $\Delta_t\subset \Theta_t$ then vertical bubbles do not appear at the limit. 

Indeed if there were such a bubble over a point $z\in S^1$ it would fill out a bounded component of the complement of $N_z$ in the fiber. In particular its interior would intersect $C$. We would still have $\Delta_t\cap C\neq \emptyset$ for $t$ close to 1 by persistence of intersection (see \S 1 b)). But this is impossible as $\Delta_t\subset \partial \hat N_t$ and $C\vert_D$ sits in the interior of $\hat N_t$. Therefore passing to the limit we get a filling $\Theta$ of $N$ made out of 1-sheeted discs $\Delta$ of index 1. 

\null\noindent
{\bf Boundary of the filling.}
We analyze now the boundaries $\partial \Delta$ on $N$. 

For this we parametrize $N$ by the standard torus $T=S^1\times S^1$ via an immersion (which is 2 to 1 precisely where $N$ is pinched). By construction the boundaries $\partial \Delta$ lift on $T$ as a foliation. We also lift $\gamma$ by doubling it.

 Homologically $\partial \Delta$ cuts $\gamma$ once in $T$. This can be checked in the model situation where the $\partial \Delta$ correspond to the constant discs $(w=u), u\in S^1$ and $\gamma$ to $(w^2=z)$ (as indices and degrees for $\pi$ agree). The point is that $\partial \Delta$ and $\gamma$ intersect exactly once {\it geometrically}. 

If not as the $\partial \Delta$ foliate $T$ we would get a contact between one of the $\partial \Delta$ and $\gamma$. Looking back in $N$ this would produce an interior intersection between $\Delta_0$ and a neighbor of $\Delta$ by persistence (see \S 1 b)). It would even persist between $\Delta_0$ and some $\Delta_t$ for $t$ close to 1. But this is impossible as $\Delta_t\subset \partial\hat N_t$ and $\Delta_0$ sits in the interior of $\hat N_t$. 

Therefore through each point of $\gamma$ pass exactly 2 boundaries $\partial \Delta$, one for each side of $T$ pinched in $N$. Note that they cross $\gamma$ in opposite directions as the sides of $T$ come with opposite orientations in the pinching.

\null\noindent
{\bf Surgery.} We focus on the point $1$ of $\gamma$. Note that $B_1 \cap (\vert w \vert \leq 1)=[-1,1]$. Denote by $\Delta_\pm$ the 2 discs of $\Theta$ whose boundaries cross at $1$. Precisely $\partial \Delta_+$ comes from above $[-1,1]$ and $\partial \Delta_-$ from below (before the pinching). These discs are of the form $(z, f_\pm(z))$ where $f_\pm$ is holomorphic in $D$ and smooth up to the boundary. We have $f_\pm(1)=0$ and $f_\pm$ does not vanish elsewhere. Moreover $if'_+(1)<0$ and $if'_-(1)>0$.

We put $\Delta'=((w-f_+(z))(w-f_-(z))=\epsilon z)$ for $\epsilon>0$ small. Over $D$ it is an embedded 2-sheeted disc close to $\Delta_+\cup\Delta_-$. Moreover $\partial \Delta'$ is still contained in the standard M\" obius strip near $1$ and avoids $\gamma$. It locally looks like a hyperbola obtained by deforming the cross $\partial \Delta_+\cup \partial \Delta_-$.

\null\noindent
{\bf Constructing $M'$.} Remove a small neighborhood of $\gamma$ in $Z$. Then thicken slightly what remains and take the boundary. We get a 2-sheeted torus $M''$. Each component of $M''_z$ looks like an enlarged component of $M_z$ to which we grafted a thin needle (with rounded tip) along the main part of half of the arc $B_z$. The torus $M'$ will be obtained from $M''$ by slightly deforming it in order to contain $\partial \Delta'$. 

This certainly can be done far from $1$. Near $1$ we deal with the ends of the needles. At $1$ they are horizontal, left and right from $0$. Let us take a look at the right branch of the hyperbola of $\partial \Delta'$. Before $1$ (for $z=e^{i\theta}$ with small $\theta<0$) we insert the branch in the bottom boundary of the right needle (we push it slightly up). At $1$ we pass its rounded tip. After $1$ ((for $z=e^{i\theta}$ with small $\theta>0$) we insert the branch in the top boundary of the needle (we push it slightly down). 

In particular when we pass the tip we get a negative half turn for $\mu(\Delta')$ with respect to $M'$ compared to $\mu(\Delta_\pm)$ with respect to $N$. The same holds for the left branch of the hyperbola. Therefore $\mu(\Delta')=\mu(\Delta_+)+\mu(\Delta_-)-1=1$. Hence $\Delta'$ generates a filling $\Sigma'$ of $M'$ (see \S 1 c)) .

\null\noindent
{\bf Deforming back to $M$.} Construct a decreasing deformation $(M_t)_{t\in [0,1]}$ from $M_0=M'$ to $M_1=M$. This can be achieved by resorbing progressively the needles in each fiber. Denote by $Y_t$ the corresponding vertical fillings. We want to produce a global deformation of (horizontal) fillings starting from $\Sigma_0=\Sigma'$. Note that $\hat Y_t\vert_D$ is non degenerate for $t<1$ as $Y$ sits in the interior of $Y_t$.

Note also that for some definite time (non depending on $\epsilon$), say $t< 2\delta$, the algebraic curve $C$ intersects $Y_t$ in its interior. As above this forbids vertical bubbles and gives compactness for fillings if $t<2 \delta$. So by the existence of local deformations (see \S 1 c)) and this compactness we get fillings $\Sigma_t$ of $M_t$ up to $t<2\delta$. 

Next we claim that $\Sigma_\delta=\partial \hat Y_\delta\vert_D$. This follows from Oka principle (see \S 1 a)). Indeed consider a disc $\Delta_t\subset \Sigma_{t}$ for $t<\delta$, follow it through the previous fillings up to a disc of $\Sigma'$. The latter deforms inside $\Sigma'$ to $\Delta'$ which deforms to $\Delta_+\cup \Delta_-$ by letting $\epsilon \to 0$. Follow this pair through the fillings $\Theta_t$ up to a pair of constant discs contained in $(\vert w\vert=R)$. They clearly avoid $\hat Y_\delta$. As the boundaries of all these discs avoid $M_\delta$ (for $\epsilon$ small) we get that $\Sigma_ {t}\cap \hat Y_\delta=\emptyset$. This concludes as the $\Sigma_t$ foliate a 1-sided neighborhood of $\Sigma_\delta$. 

Starting now from $\Sigma_\delta$ we are in position to use the compactness result for fillings as long as $\hat Y_t$ is non degenerate (see \S 1 c)). Therefore we extend the deformation of fillings up to $t<1$ by compactness and local deformations. Passing to the limit either we get a filling $ \Sigma=\partial \hat Y\vert_D$ (absence of bubbles) or $\hat Y\vert_D$ is degenerate (presence of bubbles).

\null\noindent
b) {\bf Rest of the hull.} 

We explain now the presence of a horizontal disc $\Delta$ with $\partial \Delta \subset M$ passing through a given point $x_0=(z_0,w_0)\in \hat Y\vert_D$.  Note that by the previous paragraph we only have to deal with the non degenerate case and we may suppose $x_0$ in the interior of $\hat Y\vert_D$. We proceed again by continuity. We construct a 2-sheeted torus $M'$ for which we know the presence of a disc $\Delta'$ passing through $x_0$ and deform it back to $M$.

\null\noindent
{\bf Constructing $M'$.} 
Take a 2-sheeted disc $\Delta'=((w-w_0)^2=R(z-z_0))$ through $x_0$ with $R$ large so that $\partial \Delta'$ avoids $Y$ by far.
Glue an annulus $A$ to $Y\cup\partial \Delta'$. Its fiber $A_z$ consists in 2 arcs, each connecting a point of $\partial \Delta'_z$ to a component of $Y_z$. This is possible because $w(M)=1$. Then thicken slightly $Y\cup A$ and take the boundary to get a new 2-sheeted torus $M'$ with $\partial \Delta'\subset M'$. Each component of $M'_z$ looks like an enlarged component of $M_z$ to which we grafted a thin needle reaching out the corresponding point of $\partial \Delta'_z$. Note that the index of $\Delta'$ with respect to $M'$ is $3$.  

\null\noindent
{\bf Deforming back to $M$.} Construct a decreasing deformation $(M_t)_{t\in [0,1]}$ from $M_0=M'$ to $M_1=M$. This can be achieved by resorbing progressively the needles in each fiber. Denote by $Y_t$ the corresponding vertical fillings. Note that $\hat Y_t\vert_D$ is non degenerate for all $t$. So we get (horizontal) fillings $\Sigma_t=\partial\hat Y_t\vert_D$ for all $t$ (see \S 2 a)).

 We want to produce a deformation of $\Delta'$ passing through $x_0$. This can be done locally as the index of $\Delta'$ is $>1$ (see \S 1 b)). In general denote by $I$ the set of $t$ in $[0,1]$ such that there exists a horizontal disc $\Delta_t$ passing through $x_0$ with $\partial \Delta_t\subset M_t$ and $\mu(\Delta_t)=2$ or $3$. Then $I$ is open non empty. It remains to prove that it is closed to conclude.

\null\noindent
{\bf Compactness.} Consider a sequence of such discs $\Delta_{t_n}$, $t_n\in I$, $t_n\to t$. After extracting pass to the limit thanks to Gromov theorem (see \S 1 b)). We get a horizontal limit $\Delta_t$ with the same features except maybe for its index $\mu_t$. By construction we have $\mu_t\leq 3$. Now if $\mu_t<1$ then $\Delta_t$ would be the only horizontal disc attached to $M_t$ by uniqueness (see \S 1 b)). But this is impossible as we have already a lot of other such discs, those of $\Sigma_t$. If $\mu_t=1$ then $\Delta_t$ would be part of $\Sigma_t$ by uniqueness of the filling (see \S 1 c)). But this is impossible as $\Sigma_t=\partial \hat Y_t\vert_D$ and $x_0$ sits in the interior of $\hat Y_t$. So $\mu_t=2$ or $3$ and $I$ is closed. 

\Refs 

\widestnumber\no{99}
\refno=0

\bref \by M. Audin and J. Lafontaine (ed.)\book Holomorphic curves in symplectic geometry \bookinfo Prog. Math. \vol117 \publ Birkh\"auser \yr1994 \publaddr Basel
\endref

\bref \by E. Bedford\paper Stability of the polynomial hull of $T^2$\jour Ann. Scuola Norm. Sup. Pisa\vol8\yr1981\pages311--315
\endref

\bref \by E. Bedford and W. Klingenberg\paper On the envelope of holomorphy of a 2-sphere in $\C^2$\jour J. Amer. Math. Soc.\vol4\yr1991\pages623--646
\endref

\bref \by E. M. Chirka\paper Regularity of the boundaries of analytic sets\jour Mat. Sb. vol117\yr1982\pages291--336
\endref

\bref \by J. Duval\paper Une contrainte g\'eom\'etrique pour certaines sous-vari\'et\'es rationnellement convexes\jour Math. Ann. \vol 289 \yr 1991\pages627--629
\endref

\bref \by J. Duval and N. Sibony \paper Hulls and positive closed currents \jour Duke Math. J. \vol95\yr1998\pages621--633
\endref

\bref \by Y. Eliashberg\book Filling by holomorphic discs and its applications, {\rm in} Geometry of low dimensional manifolds
\pages45--67 \bookinfo L.M.S. Lecture Notes\vol151 \publ Cambridge U.P. \yr1990 \publaddr Cambridge
\endref

\bref \by F. Forstneric\paper Analytic disks with boundaries in a maximal real submanifold of $\C^2$\jour Ann. Inst. Fourier \vol37\yr 1987\pages1--44
\endref

\bref \by F. Forstneric\paper Polynomial hulls of sets fibered over the circle\jour Indiana Univ. Math. J.\vol37\yr1988\pages869--889
\endref

\bref \by M. Gromov\paper Pseudoholomorphic curves in symplectic manifolds\jour Invent. Math.\vol82\yr1985\pages307--347
\endref

\bref\ \by S. Ivashkovich and V. Shevchishin\paper Reflection principle and $J$-complex curves with boundary on totally real immersions\jour Commun. Contemp. Math.\vol4\yr2002\pages65--106
\endref

\bref \by Z. Slodkowski\paper Polynomial hulls in $\C^2$ and quasicircles \jour Ann. Scuola Norm. Sup. Pisa \vol16\yr1990\pages367--391
\endref

\bref \by E. L. Stout \book Polynomial convexity \bookinfo Prog. Math. \vol261 \publ Birkh\"auser \yr2007 \publaddr Boston
\endref

\bref \by M. A. Whittlesey\paper Riemann surfaces in fibered polynomial hulls \jour Ark. Mat.\vol37\yr1999\pages409--423
\endref

\bref \by R. Ye\paper Filling by holomorphic curves in symplectic 4-manifolds \jour Trans. Amer. Math. Soc. \vol350\yr1998\pages213--250
\endref 
 
\endRefs

\enddocument